\documentclass[reqno,centertags,12pt]{amsart}
\usepackage{amsmath,amsthm,amscd,amssymb,latexsym,verbatim}

\usepackage{graphicx,epsf,cite}

-\marginparwidth \oddsidemargin -4mm \evensidemargin \oddsidemargin

\newtheorem{theorem}{Theorem}[section]

\theoremstyle{definition}
\newtheorem{definition}[theorem]{Definition}

\theoremstyle{remark}

\newcounter{smalllist}

\allowdisplaybreaks
\numberwithin{equation}{section}

\newcommand{\lb}{\label}

\newcommand{\beq}{\begin{equation}}
\newcommand{\eeq}{\end{equation}}

\newcommand{\bal}{\begin{align}}
\newcommand{\eal}{\end{align}}
\newcommand{\bals}{\begin{align*}}
\newcommand{\eals}{\end{align*}}

\newcommand{\bbN}{{\mathbb{N}}}
\newcommand{\bbR}{{\mathbb{R}}}

\newcommand{\bbP}{{\mathbb{P}}}

\newcommand{\bbZ}{{\mathbb{Z}}}

\newcommand{\bbQ}{{\mathbb{Q}}}

\newcommand{\bbS}{{\mathbb{S}}}

\newcommand{\calS}{{\mathcal S}}

\newcommand{\calL}{{\mathcal L}}
\newcommand{\calK}{{\mathcal K}}

\newcommand{\calF}{{\mathcal F}}

\newcommand{\eps}{\varepsilon}

\begin{document}
\title[Homogenization for KPP Reactions]
{Homogenization for Time-periodic KPP Reactions}

\author{Andrej Zlato\v s}
\address{\noindent Department of Mathematics \\ UC San Diego \\ La Jolla, CA 92093, USA \newline Email:
zlatos@ucsd.edu}


\begin{abstract}
We prove homogenization for reaction-advection-diffusion equations with KPP reactions, in the  time-periodic spatially stationary ergodic setting, and find  an explicit formula for the homogenized dynamic.  We also extend this result to models with non-local diffusion and KPP reactions. 
\end{abstract}

\maketitle

\section{Introduction and Main Results} \lb{S1}

Many  processes  occurring in random media are modeled by the reaction-diffusion equation
\beq\lb{1.4}
u_t = \calL_\omega u + f(t,x,u,\omega).
\eeq
Here $\omega$ is an element from some probability space $(\Omega, \bbP,\calF)$, and
\beq\lb{1.2b}
\calL_\omega u(t,x):= \sum_{i,j=1}^d A_{ij}(t,x,\omega) u_{x_ix_j}(t,x)  + \sum_{i=1}^d b_i(t,x,\omega) u_{x_i} (t,x)
\eeq
is a second-order differential operator, with the spatial Laplacian $\calL_\omega=\Delta$ being the basic (deterministic) case.
The unknown function $u$ can represent concentration of a substance or density of a species, which is subject to  diffusion and advection (modeled by $\calL_\omega$) as well as some reactive process (modeled by $f$), both of which can be space-time dependent.
We will study  this model and its non-local version
\beq\lb{1.2c}
\calL_\omega u (t,x):=  {\rm p.v.}  \int_{\bbR^d}  K (t,x,\nu,\omega) [u(t,x+\nu)-u(t,x)] d\nu
\eeq
with KPP (a.k.a.~Fisher-KPP)  reactions $f$.  Named after Kolmogorov, Petrovskii, and Piskunov \cite{KPP} and Fisher \cite{Fisher}, who first studied them in 1937, these  are defined as follows.

\begin{definition} \lb{D.1.0}
A Lipschitz-in-$(t,x,u)$ function $f:\bbR^+\times\bbR^d\times [0,1]\times\Omega\to\bbR$ is a {\it KPP reaction} if $f(\cdot,\cdot,0,\cdot)\equiv 0\equiv f(\cdot,\cdot,1,\cdot)$ and $f(\cdot,\cdot,u,\cdot)\le f_u(\cdot,\cdot,0,\cdot)u $ for all $u\in [0,1]$ (with $f_u(\cdot,\cdot,0,\cdot)$ existing pointwise), plus the following uniform hypotheses are satisfied. 
We have $\inf_{(t,x,\omega)\in\bbR^+\times\bbR^d\times\Omega} f(t,x,u,\omega)>0$ for each $u\in(0,1)$, plus $\inf_{(t,x,\omega)\in\bbR^+\times\bbR^d\times\Omega} f_u(t,x,0,\omega)>0$ and   
\[
\lim_{u\to 0} \sup_{(t,x,\omega)\in\bbR^+\times\bbR^d\times \Omega} \left( f_u(t,x,0,\omega)-\frac{f(t,x,u,\omega)}u \right)=0.
\]
\end{definition}

We will consider here the Cauchy problem for \eqref{1.4} with $(t,x)\in \bbR^+\times\bbR^d$, and $f(t,x,\cdot)$ will only be defined on $[0,1]$ because we are interested in solutions with normalized concentration $0\le u\le 1$.  While both $u\equiv 0$ and $u\equiv 1$ are stationary solutions for \eqref{1.4}, the former is unstable and the latter is asymptotically stable because $f>0$ when $u\in (0,1)$.  We will be interested in the transition of solutions from values near 0 to those near 1, which models invasion of the low concentration region $u\approx 0$ by the modeled substance or population.  

We will prove that {\it homogenization} --- existence of a deterministic homogenous large scale solution dynamic --  occurs for \eqref{1.4} with KPP reactions under suitable hypotheses.
Specifically, we will consider uniformly elliptic diffusion matrices $A$ and advection vectors $b$ that are bounded above by twice the square root of the product of the ellipticity constant of $A$ and $\inf f_u(\cdot,\cdot,0,\cdot)$, which will allow us to apply the concept of virtual linearity for KPP reaction-diffusion equations from the companion paper \cite{ZlaKPPlinearity}.  The latter shows that leading order of the solution dynamic for \eqref{1.4}, which is what matters for homogenization, only depends on $f$ through $f_u(\cdot,\cdot,0,\cdot)$.
We will also assume that $(A,b,f_u(\cdot,\cdot,0,\cdot))$ is time-periodic and {\it spatially stationary ergodic}.
The latter means that $(\Omega, \bbP,\calF)$ is endowed with a group of measure-preserving bijections $\{{\Upsilon_y:\Omega\to\Omega}\}_{y\in\bbR^d}$ such that $\Upsilon_0={\rm Id}_{\Omega}$,  for any $y,z\in \bbR^d$ we have
\[
\Upsilon_y\circ\Upsilon_z=\Upsilon_{y+z},
\]
and for each $(t,x,y,u,\omega)\in \bbR^+\times \bbR^{2d}\times [0,1]\times\Omega$ we have
\[
(A,b,f_u(\cdot,\cdot,0,\cdot))(t,x,\Upsilon_y\omega) = (A,b,f_u(\cdot,\cdot,0,\cdot))(t,x+y,\omega)
\] 
(this is stationarity), 
as well as that $\bbP[E]\in\{0,1\}$ holds for each $E\in\calF$ such that $\Upsilon_y E=E$ for all $y\in\bbR^d$.  
In general, stationary ergodicity is the minimal hypothesis under which one can hope for  homogenization to hold, although our results show that in the case of KPP reactions, $f$ itself will not need to be stationary ergodic.  

Since \eqref{1.4} exhibits ballistic propagation of solutions (i.e., linear rate of invasion of the region $u\approx 0$ by the region $u\approx 1$), homogenization for it involves the ballistic scaling
\begin{equation}\label{1.5}
u^{\eps}(t,x, \omega):=u\left( \eps^{-1}  t, \eps^{-1}  x, \omega\right)
\end{equation}
with $\eps>0$ small.  This turns \eqref{1.4} into its large-space-time-scale version
\begin{equation}\label{1.6}
u^{\eps}_{t}= \calL_\omega^\eps u^{\eps}+ \eps^{-1} f\left(\eps^{-1} t,\eps^{-1} x, u^{\eps}, \omega\right),
\end{equation}
where
\[
\calL_\omega^\eps u^\eps(t,x):= \eps \sum_{i,j=1}^d A_{ij} \left(\eps^{-1} t,\eps^{-1} x,\omega \right) u^\eps_{x_ix_j}(t,x)  + \sum_{i=1}^d b_i \left(\eps^{-1} t,\eps^{-1} x,\omega \right) u^\eps_{x_i} (t,x).
\]
One then hopes  that for almost every $\omega\in\Omega$, solutions to \eqref{1.6} with some $\eps$-independent initial value $\bar u(0,\cdot)$ converge as $\eps\to 0$  to some function $\bar u$ that solves some homogeneous PDE with the same initial value.  

We will show that this is indeed the case for KPP reactions, and  that the deterministic dynamic is in fact quite simple.   Essentially, if $G$ is the support of $\bar u(0,\cdot)$, then $\bar u(t,\cdot)=\chi_{G+t\calS}$ for some bounded convex open set $\calS$ that only depends on $(A,b,f_u(\cdot,\cdot,0,\cdot))$.  We note that this means (see, e.g., \cite{LinZla}) that this $\bar u$ is also a discontinuous viscosity solution to the homogeneous  Hamilton-Jacobi equation
\beq\lb{1.9a}
\bar u_t=c^* \left(-\nabla\bar u |\nabla \bar u|^{-1} \right) |\nabla \bar u|
\eeq
with
\beq\lb{1.9}
c^{*}(e):=\sup_{y\in \calS} \,y\cdot e.
\eeq
The set $\calS$ 
can be obtained as an asymptotic growth shape of level sets of initially compactly supported solutions to \eqref{1.4}, per the following definition whose second part is from \cite{LinZla}.

\begin{definition} \label{D.1.3}
If for some $\omega\in \Omega$ there is an open set $\mathcal S\subseteq\bbR^d$ such that the solution 
$u(\cdot,\cdot,\omega)$ 
to \eqref{1.4} with $u(0,\cdot,\omega):=\frac 12\chi_{B_1(0)}$ satisfies
\begin{align}
\lim_{t\rightarrow\infty}  
\inf_{x\in (1-\delta) \mathcal{S}t} u(t,x,\omega) & =1, \label{1.3}
\\ \lim_{t\rightarrow\infty} 
\sup_{x\notin (1+\delta) \mathcal{S}t} u(t,x,\omega) & =0 \label{1.3a}
\end{align}
for  each $\delta>0$, 
then  $\mathcal S$ is a {\it Wulff shape} for \eqref{1.4} with this $\omega$. If there is $\Omega'\subseteq \Omega$ with $\bbP[\Omega']=1$ such that \eqref{1.4} with each $\omega\in\Omega'$ has the same Wulff shape $\mathcal S$, then  $\mathcal S$ is a {\it deterministic Wulff shape} for \eqref{1.4}.

Similarly, if for some $(\omega,e)\in \Omega \times \bbS^{d-1}$ there is $c^*(e)$ such that the solution $u(\cdot,\cdot,\omega;e)$  to \eqref{1.4} with  $u(0,\cdot,\omega;e):=\frac 12\chi_{\{x\cdot e<0\}}$ satisfies
\begin{align}
\lim_{t\rightarrow\infty}  
\inf_{x\in (c^*(e)e-Q)t} u(t,x,\omega;e) & =1, \label{1.3b}
\\ \lim_{t\rightarrow\infty} 
\sup_{x\in (c^*(e)e+Q)t}  u(t,x,\omega;e) & =0 \label{1.3c}
\end{align}
for each compact set $Q\subseteq \{x\cdot e>0\}$, then   $c^*(e)$ is a {\it front speed} in direction $e$ for \eqref{1.4} with this $\omega$. If there is $\Omega'\subseteq \Omega$ with $\bbP[\Omega']=1$ such that \eqref{1.4} with each $\omega\in\Omega'$ has the same front speed $c^*(e)$ in direction $e$, then  $c^*(e)$ is a {\it deterministic front speed} in direction $e$ for \eqref{1.4}.
\end{definition}

The simplicity of the formula $\bar u(t,\cdot)=\chi_{G+t\calS}$ for the homogenized dynamic is due to the second main contribution of \cite{ZlaKPPlinearity}, which shows that solutions to \eqref{1.4} with KPP reactions exhibit yet another form of linear-like behavior.  Namely,  the leading order of the solution dynamic for a general initial datum $u(0,\cdot)$ can be recovered from solving the PDE with initial data that are obtained by restricting $u(0,\cdot)$ to the cubes $C_n:= (n_1,n_1+1)\times\dots\times(n_d,n_d+1)$ with $n\in \bbZ^d$, so the nonlinear interaction between the resulting initially compactly supported solutions does not affect the leading order of the solution dynamic.  This mean that we will be able to recover the above formula for general initial data from existence of a Wulff shape for \eqref{1.4}.  For the reader's convenience, we provide here the random medium version of the main virtual linearity result from \cite{ZlaKPPlinearity} (Theorem 1.2 there; see also Remark 2 following it).

\begin{theorem} \lb{T.1.1}
Let  $f$ be a KPP reaction 
and $\calL_\omega$ be from \eqref{1.2b}, where 
 $A=(A_{ij})$ is a bounded symmetric matrix with $A\ge \lambda I$ for some $\lambda>0$ (and each $\omega\in\Omega$), and the vector $b=(b_1,\dots,b_d)$ satisfies $\|b\|_{L^\infty}^2< 4\lambda \inf_{(t,x,\omega)\in\bbR^+\times\bbR^d\times\Omega} f_u(t,x,0,\omega)$.
Let 
 \beq \lb{2.30}
 f'(t,x,u,\omega):=f_u(t,x,0,\omega)\min\{u,1-u\}
 \eeq
  for $(t,x,u,\omega)\in \bbR^+\times\bbR^d\times[0,1]\times\Omega$.
Then there is $\phi :\bbR^+\to\bbR^+$ with $\lim_{s\to \infty} \phi (s)=0$, 
and for each  $\delta\in(0,\frac 12]$ there is $\tau_\delta\ge 1$, 
such that  the following holds for any fixed $\omega\in\Omega$.

If  $u:\bbR^+\times \bbR^d\to [0,1]$  solves \eqref{1.4}, and for each $n\in\bbZ^{d}$ we let $u_n' :\bbR^+\times \bbR^d\to[0,1]$ solve \eqref{1.4} with $f'$ in place of $f$ and with $u_n' (0,\cdot):=  u(0,\cdot) \chi_{C_n}$, 
then for each $(t,x)\in[\tau_\delta,\infty)\times \bbR^d$,
\beq  \lb{2.11}
\sup_{n\in\bbZ^d} u_n' \left( t-\delta t, x \right) -\phi(\delta t) \le u(t,x)  \le  \sup_{n\in\bbZ^d}  u_n'  \left(t+\delta t,x \right)  +\phi(\delta t).
\eeq
\end{theorem}

{\it Remark.} 
So up to $o(t)$ time shifts  and $o(1)$ errors, we have  $u\approx \sup_{n\in\bbZ^d} u_n' $. 
\smallskip

The following is our main homogenization result, in which we allow the initial values to also have $\eps$-dependent shifts $y_\eps$ plus perturbations whose sizes decay as $\eps\to 0$.  We let here
 $B_r(G):=G+(B_r(0)\cup\{0\})$ and $G^0_r:=G\backslash\overline{B_r(\partial G)}$ for any $G\subseteq\bbR^d$ and $r\ge 0$ (so $G^0_0$ is the interior of $G$).

\begin{theorem} \lb{T.1.2}
Assume that the hypotheses of Theorem \ref{T.1.1} hold, and $(A,b,f_u(\cdot,\cdot,0,\cdot))$ is 
 time-periodic and  spatially stationary ergodic.
 Then \eqref{1.4} with $f$ as well as with $f'$ from \eqref{2.30}
 has the same deterministic convex bounded  Wulff shape $\calS\ni 0$, and for almost all $\omega\in\Omega$ the following holds.
 If  $G\subseteq\bbR^d$ is open, $\theta\in(0,1)$, $\Lambda>0$, and $u^\eps(\cdot,\cdot,\omega)$ solves \eqref{1.6} with
\begin{equation} \label{1.7}
\theta\chi_{(G+y_\eps)_{\rho(\eps)}^0} \le u^\eps(0,\cdot,\omega)\le \chi_{B_{\rho(\eps)}(G+y_\eps)}
\end{equation}
for each  $\eps>0$, with some $y_\eps\in B_{\Lambda}(0)$ and $\lim_{\eps\to 0}\rho(\eps)=0$ (when $y_\eps=0$ and $\rho(\eps)= 0$, this becomes just $\theta\chi_{G} \le u^\eps(0,\cdot,\omega)\le \chi_{G}$), then 
\begin{equation} \label {1.8} 
\lim_{\eps\to 0}  u^\eps(t,x+y_\eps,\omega)= \chi_{G^\calS}(t,x)
\end{equation}
locally uniformly on $ ([0,\infty)\times\bbR^d) \setminus \partial G^\calS$, where $G^\calS:=\{(t,x)\in\bbR^+\times\bbR^d \,|\, x\in G+t\calS\}$.
In particular \eqref{1.4} has a deterministic front speed $c^*(e)$ for each $e\in\bbS^{d-1}$, given by \eqref{1.9}.
\end{theorem}

{\it Remark.}  
%
Theorems 2.6 and A.2 in the paper \cite{CLM} by Caffarelli, Lee, and Mellet imply existence of a time-independent space-periodic (so non-random) ignition reaction in two dimensions  and of $e\in\bbS^1$ such that $\calS$ and $c^*(e)$ as above exist  when $\calL_\omega=\Delta$, but $c^{*}(e)>\sup_{y\in \calS} \,y\cdot e$.  This shows that Theorem \ref{T.1.2} does not hold for general non-KPP reactions $f$, and neither does the second inequality in \eqref{2.11} (even if we let $f':=f$).  Nevertheless, homogenization results for ignition reactions in $\bbR^d$ for $d\le 3$ were obtained by Lin and the author  \cite{LinZla}, and by Zhang and the author \cite{ZhaZla, ZhaZla3}.\smallskip


We note that homogenization for KPP reactions with {\it time-independent} stationary ergodic $(A,b,f)$ is stated in Theorem 9.3 of the paper \cite{LioSou} by Lions and Souganidis; however, while they indicate that a proof can be obtained via methods developed in \cite{LioSou} and three other papers, it is not provided there.  That result is stated with the right-hand side of \eqref{1.8} being the characteristic function of  ${\rm int}\, Z^{-1}(0)$, where $Z\ge 0$ is the viscosity solution to $\min\{Z_t+H(\nabla Z),Z\}=0$ on $\bbR^+\times\bbR^d$ with initial value $Z(0,\cdot) = \infty\chi_{\bbR^d\setminus G}$, and $H$ is some (non-explicit) homogeneous Hamiltonian depending on $(A,b,f_u(\cdot,0,\cdot))$.
While the Wulff shape $\calS$ in our Theorem \ref{T.1.2} also depends on $(A,b,f_u(\cdot,\cdot,0,\cdot))$, \eqref{1.8} yields a simpler and more explicit description of the limiting solution dynamic  (also, instead of the above PDE for $Z$, we now have \eqref{1.9a} with a 1-homogeneous Hamiltonian given by \eqref{1.9}).  Moreover, our approach applies in the time-periodic $(A,b,f_u(\cdot,\cdot,0,\cdot))$ case, and also extends to \eqref{1.4} with time-periodic non-local diffusion operators (see Theorem \ref{T.1.1a} below).

Besides \cite[Theorem 9.3]{LioSou}, we are not aware of other prior results on homogenization for KPP reactions in dimensions $d\ge 2$, even in the case of periodic reactions (although existence of Wulff shapes and front speeds in the periodic case goes back to \cite{GF}).   However, multidimensional homogenization results have been proved recently by the author and others for other types of reactions:  for periodic ignition and monostable reactions in  \cite{LinZla, AlfGil}, and for general stationary ergodic ignition reactions in \cite{LinZla, ZhaZla, ZhaZla3} (see also \cite{BarSou, MajSou} for the cases of homogeneous KPP or bistable reactions and periodic first- or second-order terms).  All these results except for \cite{MajSou} also concern the time-independent coefficients/reactions case.

We also add that Zhang and the author were able to apply Theorem \ref{T.1.1} to the problem of homogenization in certain environments with general space-time dependence.  Specifically, the companion paper \cite{ZhaZla5} proves Theorem \ref{T.1.2} for $(A,b,f_u(\cdot,\cdot,0,\cdot))$ 
that has temporally decaying correlations in an appropriate sense (this is essentially the opposite of the  time-periodic case from Theorem \ref{T.1.2}, which is fully time-correlated).

\medskip


We also show that Theorem \ref{T.1.2} extends to \eqref{1.4} with non-local diffusion operators from \eqref{1.2c} under suitable hypotheses  (see \cite{ZlaKPPlinearity}).   In particular, for the solution dynamic to be well-defined and ballistic in the presence of KPP reactions, the diffusion kernel $K$ must decay at least exponentially as $\nu\to\infty$ and have $O(|\nu|^{-d-2+\alpha})$ growth as $\nu\to 0$ (for some $\alpha>0$).  This growth coincides with $\calL_\omega=-(-\Delta)^{\alpha/2}$, for which well-posedness, comparison principle, and the parabolic Harnack inequality are known to hold \cite{BSV}.  The methods used to establish these properties should equally apply to general kernels with the above asymptotics, and we will assume them as hypotheses below so that our result applies whenever these can be established.  We note that the Harnack inequality referred to here is the forward one (equations with non-local diffusions may also satisfy backward-in-time Harnack inequalities, such as in \cite{BSV}).

For kernels as above, scaling \eqref{1.5} leads to the operators
\beq\lb{1.2e}
\calL_\omega^\eps u^\eps(t,x):= \eps^{-d-1} \, {\rm p.v.} \int_{\bbR^d}  K  \left(\eps^{-1}t,\eps^{-1}x,\eps^{-1}\nu,\omega \right) [u^\eps(t,x+\nu)-u^\eps(t,x)] d\nu ,
\eeq
and we then obtain the following non-local version of Theorem \ref{T.1.2}.

\begin{theorem} \lb{T.1.1a}
Let $\alpha\in(0,1]$ and let $F$ be some family of even-in-$\nu$ kernels such that for each $ K \in F$,  there is $\calK:(0,\infty)\to[0,\infty)$ with $\chi_{(0,\alpha]}(r)\le  \calK(r)\le \chi_{(0,\alpha]}(r)r^{-d-2+\alpha}$ on $(0,\infty)$ and 
\beq\lb{2.31}
\alpha  \calK(|\nu|) \le   K (t,x,\nu,\omega) \le \alpha^{-1} \max\left\{ \calK(|\nu|),  e^{-\alpha|\nu|} \right\} 
\eeq
for  each $(t,x,\nu,\omega)\in\bbR^+\times\bbR^{2d}\times\Omega$.  Assume also that \eqref{1.4} with any $ K\in F $ and $\calL_\omega$ from \eqref{1.2c}, any KPP reaction $f$, and locally BV initial data $0\le u(0,\cdot)\le 1$ is well-posed in some subspace $\mathcal A\subseteq L^{1}_{\rm loc}(\bbR^+\times \bbR^d)$, where the comparison principle for sub- and supersolutions to \eqref{1.4} as well as the parabolic (forward) Harnack inequality also hold (the latter with  uniform constants for all $K\in F$ and all $f$ with the same  Lipschitz constant), and the solutions for $u(0,\cdot)\equiv 0,1$ are $u\equiv 0,1$, respectively.   

Let $f$ be a KPP reaction, let $K\in F$, and assume that $( K ,f_u(\cdot,\cdot,0,\cdot))$ is time-periodic and  spatially (in $x$) stationary ergodic.  Then the claims in Theorem~\ref{T.1.2} hold for $\calL_\omega$ from \eqref{1.2c} and $\calL_\omega^\eps$ from \eqref{1.2e}.
\end{theorem}

{\it Remark.}
One could also extend this to mixed diffusion operators, with $\calL_\omega$ being the sum of the right-hand sides of \eqref{1.2b} and  \eqref{1.2c}, but we will not do so here.
\smallskip


\subsection{Acknowledgements}
The author thanks Scott Armstrong, Jessica Lin, James Norris,  and Yuming Paul Zhang for useful discussions and pointers to literature.  He also acknowledges partial support by  NSF grant DMS-1900943 and by a Simons Fellowship.

\section{Proofs of Theorems \ref{T.1.2} and \ref{T.1.1a}} \lb{S4}

Let us now prove Theorem \ref{T.1.2}.  This proof also applies to Theorem \ref{T.1.1a}, if we replace $(A,b)$ by $K$ and Theorem~\ref{T.1.1} by  \cite[Theorem 1.3]{ZlaKPPlinearity}. 

\medskip

Let us first assume that $(A,b,f)$ is also time-periodic and spatially stationary ergodic, and  assume for simplicity that its temporal period is 1 (the argument is the same for any period).
We  note that  for any fixed $\theta\in(0,1)$, any solution to \eqref{1.4} with $u(0,\cdot)\ge \theta \chi_{B_1(0)}$ converges locally uniformly on $\bbR^d$ to 1 as $t\to\infty$ provided the hypotheses of Theorem~\ref{T.1.1} are satisfied (this is called {\it hair-trigger effect}), and this convergence is uniform in all $(A,b,f)$ that satisfy these hypotheses uniformly (including a uniform positive lower bound on $4\lambda \inf_{(t,x,\omega)\in\bbR^+\times\bbR^d\times\Omega} f_u(t,x,0,\omega) - \|b\|_{L^\infty}^2$).  Indeed, this is proved by the argument in the proof of Theorem \ref{T.1.1} in \cite{ZlaKPPlinearity} (where it is Theorem 1.2) that verifies estimate (2.2) in \cite{ZlaKPPlinearity}. 
Hence under our hypotheses, this convergence is uniform in $\omega\in\Omega$.  This uniformity then also extends to any shift of the initial datum, after accounting for the corresponding shift in the solution, because shifting the medium by $y\in\bbR^d$ simply amounts to changing $\omega$ to $\Upsilon_y \omega$.  

Thus if for each $y\in\bbR^d$ we let $u(\cdot,\cdot,\omega;y)$ solve \eqref{1.4} with $u(0,\cdot,\omega;y)=\frac 12 \chi_{B_1(y)}$, and define 
\[
\tau(y, z, \omega):=\inf\left\{ t\in\bbN_0 \,\bigg|\, u(s,\cdot, \omega; y)\geq \frac 12\chi_{B_{1}(z)} \text{ for all $s\ge t$} \right\} \ge 0,
\]
then there is $C\ge 1$ such that $\tau(y,z,\omega)\le C$ whenever $|y-z|\le 1$.  Comparison principle shows that for all $x,y,z\in\bbR^d$ we have the subadditive estimate
\beq\lb{3.1}
\tau(y, z, \omega)\leq \tau(y, x, \omega)+\tau(x, z, \omega)
\eeq
because $(A,b,f)$ has temporal period 1. Then we also obtain
\beq\lb{3.2}
\tau(y,z,\omega) \le C(|y-z|+1)
\eeq
and more generally
\beq\lb{3.3}
|\tau(y,z,\omega)- \tau(y',z',\omega)| \le C(|y-y'|+|z-z'|+2),
\eeq
for all $(y,z,y',z',\omega)\in(\bbR^d)^4\times\Omega$.  Indeed, \eqref{3.2} follows from iterating \eqref{3.1} with intermediate points $y+n\frac{z-y}{|z-y|}$ for $n=1,\dots, \lfloor |z-y| \rfloor$, while \eqref{3.3} uses \eqref{3.2} and 
\[
\tau(y',z',\omega) \le \tau(y',y,\omega) + \tau(y,z,\omega) + \tau(z,z',\omega).
\]
Finally, stationarity of $(A,b,f)$ clearly also yields for all $(y,z,x,\omega)\in(\bbR^d)^3\times\Omega$,
\beq\lb{3.4}
\tau(y, z, \Upsilon_{x}\omega)=\tau(x+y, x+z, \omega).
\eeq

The above formulas and Kingman's subadditive theorem \cite{Kin} now show that for each $e\in\bbS^{d-1}$, there is 
a constant $w(e)>0$ (the {\it spreading speed} for \eqref{1.4} in direction $e$) and a set $\Omega_e\in\calF$ with $\bbP[\Omega_e]=1$ such that for each $\omega\in \Omega_e$ we have
\begin{equation} \label{3.5}
\lim_{n\rightarrow\infty} \frac{\tau(0, ne, \omega)}{n}=\frac 1{w(e)}.
\end{equation}
Moreover, there is  $c>0$ such that for all $e,e'\in \mathbb{S}^{d-1}$ we have
\begin{equation}\label{3.6}
c\leq w(e)\leq  c^{-1}
\end{equation}
and
\begin{equation}\label{3.7}
\max\left\{ \left| w(e) - w(e') \right|, \left| \frac 1{w(e)} - \frac 1{w(e')} \right| \right\} \leq c^{-1}|e-e'|.
\end{equation}
These claims are identical to those in \cite[Lemma 3.3]{LinZla} (which concerned time-independent reactions), as are their proofs, which we therefore skip.  We only note that $c=C^{-1}$ works in the first inequality of \eqref{3.6} as well as in the estimate of the second term on the left-hand side of  \eqref{3.7}, due to \eqref{3.2} and \eqref{3.3}.  The full estimate in \eqref{3.7} then clearly holds with $c=C^{-3}$, while the second inequality of \eqref{3.6} follows from  $v(t,x):=e^{at-(x-x_0)\cdot e}$ being a supersolution to \eqref{1.4} for any  $(x_0,e)\in\bbR^d\times \bbS^{d-1}$, where $a:= \gamma(1+d+d^2)$ and
\[
\gamma:=\max\left\{  \max_{i,j}\|A_{i,j}\|_{L^\infty}, \,  \max_i\|b_i\|_{L^\infty}, \, \|f_u(\cdot,\cdot,0)\|_{L^\infty} \right\} 
\]
(because $v$ ``travels'' at the constant speed $a$).

The claim that \eqref{1.4} has a deterministic Wulff shape $\calS$ now follows 
from \eqref{3.1}, \eqref{3.2}, \eqref{3.5}--\eqref{3.7}, and the hair-trigger effect.  
Indeed, we will show that
\[
\calS:=\left\{se \,|\, e\in \mathbb{S}^{d-1} \text{ and }  s\in\left[0,w(e)\right) \right\},
\]
which also implies boundedness of $\calS$ and $0\in\calS$, while convexity  is due to \eqref{3.1} and \eqref{3.5}.

Let us consider the full measure set  $\Omega':=\cap_{e\in S}\Omega_e \in\calF$ for some countable dense $S\subseteq\bbS^{d-1}$, and fix any $\omega\in\Omega'$.  Then for each $\delta\in(0,1)$, there is a finite subset $S'\subseteq S$ such that $\mathbb{S}^{d-1} \subseteq \bigcup_{e\in S'} B_{\delta}(e)$.  Let $N_\delta$ be such that
\begin{equation} \label{3.8}
\sup_{e\in S'}\sup_{n\ge N_\delta} \left| \frac n {\tau(0, ne, \omega) w(e)} - 1 \right| \le \delta
\end{equation}
and denote $\Gamma_\theta(t):=\{ x\in\bbR^d\,|\, u(t,x,\omega;0)\ge\theta\}$.  Thus for each $t\ge CN_\delta$ we obtain from this and \eqref{3.2} (which yields $\tau(0,ne,\omega)\le t$ for $n\in\{0,1,\dots,N_\delta-1\}$) that for all $e\in S'$ we have
\[
 \bigcup_{n\in\left[0,(1-\delta)w(e) t \right] \cap\bbN_0} B_1(ne) \subseteq \Gamma_{1/2}(t). 
\]
Then from \eqref{3.1}, \eqref{3.2},  \eqref{3.7}, and $\mathbb{S}^{d-1} \subseteq \bigcup_{e\in S'} B_{\delta}(e)$ we obtain for all $t\ge CN_\delta$ and some $\delta$-independent $C'>0$ that
\[
 t\calS\subseteq \Gamma_{1/2}((1+C'\delta)t).
\]
This and  the hair-trigger effect show that for each $\theta\in(0,1)$ there is $\delta$-independent $T_\theta>0$ such that 
for all $t'\ge CN_\delta(1+C'\delta)+T_\theta$  we have
\[
(1+C'\delta)^{-1} (t'-T_\theta)\calS\subseteq \Gamma_{\theta}(t').
\]
Since this holds for all $\delta>0$ and $\theta\in(0,1)$, we see that \eqref{1.3} holds for all $\omega\in\Omega'$ and $\delta>0$. 

To prove \eqref{1.3a},   note that \eqref{3.8} also implies for all $e\in S'$, $t\ge \frac {N_\delta}c$, and $n> (1+\delta)w(e)t$ ($\ge N$ by \eqref{3.6}) that
$B_1(ne) \not \subseteq \Gamma_{1/2}(s)$
for some $s\ge t$.  Again using \eqref{3.1}, \eqref{3.2},  \eqref{3.7}, and $\mathbb{S}^{d-1} \subseteq \bigcup_{e\in S'} B_{\delta}(e)$ yields that 
\[
B_1(x) \not \subseteq \Gamma_{1/2}((1-C'\delta)t)
\]
for some $\delta$-independent $C'>0$, all $t\ge \frac {N_\delta}c$, and all $x\notin  t\calS$.  Finally, the hair-trigger effect and parabolic Harnack inequality show that
for each $\theta\in(0,1)$ there is $\delta$-independent $T_\theta>0$ such that 
for all $t'\ge \frac {N_\delta}{c}(1-C'\delta)-T_\theta$  we have
\[
 \left[ \bbR^d \setminus (1-C'\delta)^{-1} (t'+T_\theta)\calS \right] \cap \Gamma_{\theta}(t')=\emptyset.
\]
Since this holds for all  $\delta>0$ and $\theta\in(0,1)$, we see that \eqref{1.3a} holds for all $\omega\in\Omega'$ and $\delta>0$. 

Hence we proved the first claim.  In fact, a stronger version of this claim holds, namely that for almost all $\omega\in\Omega$, the solutions $u(\cdot,\cdot,\omega;y)$ defined at the start of this proof satisfy
\begin{align}
\lim_{t\rightarrow\infty} \inf_{|y|\leq \Lambda t} \, \inf_{x\in (1-\delta)\mathcal{S}t}u(t,x+y, \omega; y) & =1, \label{3.9}
\\ \lim_{t\rightarrow\infty} \sup_{|y|\leq \Lambda t} \, \sup_{x\notin (1+\delta)\mathcal{S}t}u(t,x+y, \omega; y) & =0.  \label{3.9a}
\end{align}
for any $\Lambda,\delta>0$ (this was referred to in \cite{LinZla} as \eqref{1.4} having a {\it strong deterministic Wulff shape} $\calS$).
The proof of this is the same as that of the identical result \cite[Proposition 3.4]{LinZla} for time-independent reactions, using \eqref{1.3}, \eqref{1.3a}, \eqref{3.1}, \eqref{3.2},  and Egorov's and Wiener's Ergodic Theorems
(with $R_0:=1$, $\theta_0:=\frac 12$,  $c:=C^{-1}$, and $\lambda_{1-\theta_0}:=0$ in that proof), and no other properties of \eqref{1.4}.  We note that \eqref{1.3}, \eqref{1.3a}, and the two theorems are used to show that for each $\delta,\theta\in(0,1)$,  there is a full measure set $\Omega_{\delta,\theta}$ such that for each $\omega\in \Omega_{\delta,\theta}$, each $\Lambda>0$, and all large enough $t>0$ we have
\beq\lb{3.10}
(1-\delta)t\calS\subseteq \{ x\in\bbR^d\,|\, u(t,x+y,\omega;y)\ge\theta\} \subseteq  (1+\delta)t\calS
\eeq
for each $y$ in a subset of $B_{\Lambda t}(0)$ whose volume is at least $(1-\delta^{d})|B_{\Lambda t}(0)|$.  But then \eqref{3.1} and \eqref{3.2} show that \eqref{3.10} holds for all large $t$ and all $y\in B_{\Lambda t}(0)$ if we replace $\delta$ by $C'\Lambda\delta$ for some $C'>0$, so \eqref{3.9} and \eqref{3.9a} hold for all $\omega$ from the full-measure set $\bigcap_{\delta,\theta\in\bbQ\cap(0,1)} \Omega_{\delta,\theta}$.


Now consider the case when only $(A,b,f_u(\cdot,\cdot,0,\cdot))$ is time-periodic and stationary ergodic.  Then $(A,b,f_u(\cdot,\cdot,0,\cdot)\min\{u,1-u\})$ is also such, so all the above arguments apply to \eqref{1.4} with $f'$ from \eqref{2.30} in place of $f$.  Hence we  obtain existence of some strong deterministic Wulff shape $\calS$ in this setting.  But Theorem \ref{T.1.1} with $u(0,\cdot)=\frac 12 \chi_{B_1(0)}$ then shows that  $\calS$ must also be a deterministic Wulff shape for \eqref{1.4} with $f$ (by taking $t\to\infty$ for any fixed $\delta>0$, and then taking $\delta \to 0$).

Let now $\Omega'\in \calF$ with $\bbP[\Omega']=1$ be such that \eqref{3.9} and \eqref{3.9a} hold for all $\Lambda, \delta>0$, with $f'$ in place of $f$, and fix any $\omega\in\Omega'$.
Note that the ``unscaled'' versions of \eqref{1.7} and \eqref{1.8} are
\begin{equation} \label{1.7a}
\theta\chi_{(\eps^{-1} (G+y_\eps))_{\eps^{-1} \rho(\eps)}^0} \le u_\eps (0,\cdot,\omega)\le \chi_{B_{\eps^{-1} \rho(\eps)}(\eps^{-1} (G+y_\eps)) }
\end{equation}
and
\begin{equation} \label {1.8a} 
\lim_{\eps\to 0}  u_\eps(\eps^{-1} T,\eps^{-1} (x+ y_\eps),\omega)= \chi_{G^\calS}(T,x),
\end{equation}
with $u_\eps$ solving \eqref{1.4}.  Let us first consider the case of  bounded $G$, that is, we have $G\subseteq B_\Lambda(0)$ after possibly increasing $\Lambda$ from the statement of the theorem.   Fix any $T_0>0$.  

Applying Theorem \ref{T.1.1} to the initial values from \eqref{1.7a}, together with  \eqref{3.9} with $\frac{2\Lambda}{T_0}$ in place of $\Lambda$,  and with the fact that $\calS$ is the deterministic Wulff shape for \eqref{1.4} with $f'$ in place of $f$, shows that for any $\delta>0$ we have
\[
\lim_{\eps\to 0} \inf_{T\ge T_0} 
\inf_{z\in (\eps^{-1} (G+y_\eps))_{\eps^{-1} \rho(\eps)}^0 + \eps^{-1}T(1-\delta)^2 \calS } 
u_\eps(\eps^{-1}T,z,\omega)=1
\]
(here we also used the hair-trigger effect, which shows that solutions to \eqref{1.4} with initial data $\ge\theta\chi_{C_n}$ grow above $\frac 12 \chi_{B_1(n)}$ in uniform time).  But since the set under the $\inf$  contains $\eps^{-1} (G + T(1-\delta)^3 \calS +y_\eps)$ for any $T\ge T_0$ as long as $\eps>0$ is small enough, taking $\delta\to 0$ yields \eqref{1.8a} locally uniformly on $G^\calS$.  We can use a similar argument based on  \eqref{3.9a} to show that
\begin{equation} \label {1.8g} 
\lim_{\eps\to 0} \sup_{T\ge T_0} 
\sup_{z\notin B_{\eps^{-1} \rho(\eps)}(\eps^{-1} (G+y_\eps))  + \eps^{-1}T(1+\delta)^2 \calS } 
u_\eps(\eps^{-1}T,z,\omega)=0,
\end{equation}
provided that we also have $u_\eps (0,\cdot,\omega)\le \frac 12$.
We obviously obtain the same result for \eqref{1.4} with $f'$ in place of $f$.  But this means that we now obtain \eqref{1.8g} without the additional hypothesis $u_\eps (0,\cdot,\omega)\le \frac 12$ because $\frac 12 u_\eps$ is a  subsolution to \eqref{1.4} with $f'$ in place of $f$ that is initially  $\leq \frac 12$.
So after again taking $\delta\to 0$, we obtain \eqref{1.8a} locally uniformly on $((0,\infty)\times\bbR^d)\setminus \overline {G^\calS}$.  

To prove the second claim in Theorem~\ref{T.1.2} for all bounded $G$, it now suffices to extend this convergence to the union of the above set and $\{0\}\times(\bbR^d\setminus\overline G)$.
But this holds thanks to the uniform upper bound $a$ above on the spreading speed of solutions to \eqref{1.4}.
Indeed, the supersolutions $v(t,x):=e^{at-(x-x_0)\cdot e}$ mentioned above show that if $G$ is convex, then for each $x\notin \overline G$ and any solution $u_\eps$  to \eqref{1.4} with $u_\eps (0,\cdot,\omega)\le \chi_{B_{\eps^{-1} \rho(\eps)} (\eps^{-1} (G+y_\eps))}$ we have
\[
\lim_{\eps\to 0} \sup_{T\in[0, a^{-1}d(x,G)-\delta]}  u_\eps(\eps^{-1} T, \eps^{-1}(x+y_\eps),\omega) =0
\]
for any $\delta>0$.  If $G$ is not convex, we can obtain the same result by  either using Theorem~\ref{T.1.1} (since the cubes $C_n$ are convex) or the supersolutions $e^{a(t-t_0)}\sum_{i=1}^d (e^{(x-x_0)\cdot e_i} + e^{-(x-x_0)\cdot e_i})$ (with $\{e_1,\dots,e_d\}$ the standard basis in $\bbR^d$).

Finally, if $G$ is unbounded, then it obviously suffices to prove \eqref{1.8a} locally uniformly on $([0,M]\times B_M(0))\setminus\partial G^\calS$ for any $M>0$.  The last argument above (together with Theorem~\ref{T.1.1}) shows that if we only consider $(T,x)$ in this set, we can replace $G$ by $G\cap B_{(2+a^{-1})M}(0)$ because $u_\eps(0,\cdot,\omega)$ at points outside $B_{\eps^{-1}(2+a^{-1})M}(0)$  will have no effect on $u_\eps(\cdot,\cdot,\omega)$ on the set  $[0,\eps^{-1}M]\times B_{\eps^{-1}M}(0)$ in the limit $\eps\to 0$.  But since $G\cap B_{(2+a^{-1})M}(0)$ is bounded, the argument in the bounded $G$ case applies and yields the second claim in the theorem for unbounded $G$ as well. 


The third claim is its immediate consequence  by taking $G:=\{x\in\bbR^d\,|\, x\cdot e <0\} $. 


\end{document}